%% file: main.tex
\documentclass[10pt, a4paper, reqno]{amsart}
\usepackage{extra/style}

\begin{document}

    \input{files/front}
    \input{files/introduction}

    \input{files/section1}
    \input{files/section2}
    \input{files/section3}
    \input{files/section4}
    \input{files/section5}

\bibliographystyle{alpha}
\bibliography{extra/bibliography}

\end{document}

%% file: files/front.tex
\begin{abstract}
    In the setting of non-Noetherian formal schemes, we define the infinitesimal neighbourhood of a morphism to be the largest possible thickening through which the morphism factors.
    We show that formally unramified adic morphisms of formal schemes locally of formal finite presentation admit the existence of infinitesimal neighbourhoods locally given by completing with respect to an ideal.
    As a consequence, we show that such morphisms can be expressed as the composition of a closed immersion and a formally étale morphism.
    The framework herein provided should serve as a sufficiently general foundation for the development of a theory of formal stacks and algebraic foliations in arbitrary characteristic.
\end{abstract}

\maketitle

\tableofcontents

%% file: files/introduction.tex
\section*{Introduction}

Following work of McQuillan (\cite{MR4507257}), it is becoming clearer that foliations in algebraic geometry are better defined as formal algebraic stacks.
Such point of view, alongside an analogue of Frobenius theorem, will be explored in future work.
Before that, it seems necessary to discuss some adequate foundations for formal schemes. \par

Formal schemes are notoriously hard to work with.
In \cite[\S 10]{MR217083}, Grothendieck develops a satisfactory theory only for locally Noetherian formal schemes.
Unfortunately, this generality is not sufficient for our needs as the category of locally Noetherian formal schemes is not closed under fibre products.
Much later, McQuillan (\cite{MR1941631}) and Yasuda (\cite{MR2520785}) extended some results to more general formal schemes, however this requires working outside the category of locally topologically ringed spaces and is thus technically much more involved.
The approach employed in this article is to only consider formal schemes whose local models are adic rings with a \emph{finitely generated} ideal of definition.
Such formal schemes are general enough to be closed under fibre products, unlike Noetherian formal schemes, and yet restrictive enough for important properties to be preserved by base change, unlike admissible formal schemes. \par

Formal schemes naturally arise as infinitesimal neighbourhoods of morphisms of schemes.
This is the main object of interest in this article.
In differential geometry, given an embedding of smooth manifolds $f : X \hookrightarrow Y$, the tubular neighbourhood theorem asserts that the normal bundle of $f$ has an open subset containing $X$ which is diffeomorphic to an open subset $U$ of $Y$ containing $X$.
When $f$ is a morphism of algebraic varieties, it is hopeless to find such an open subset $U$ as the Zariski and \'{e}tale topologies are too coarse.
On the other hand, by \emph{completing} $Y$ along $X$, we get a factorisation $X \rightarrow \hat{Y} \rightarrow Y$, where the object $\hat{Y}$ shares similar properties with $U$.
However $\hat{Y}$ is not in general a scheme, it is a \emph{formal scheme}. \par

The morphism $X \rightarrow \hat{Y}$ is a closed immersion of formal schemes inducing a homeomorphism of topological spaces.
This is defined as a \emph{thickening} of formal schemes (see Definition \ref{def:thickening}).
We then define the \emph{infinitesimal neighbourhood} of $f$ to be the final object in the category of thickenings over $f : X \rightarrow Y$, when it exists (see Definition \ref{def:infinitesimal_neighbourhood}). \par

The main result of this article is the following.

\begin{mainthm}
    Let $f : X \rightarrow Y$ be a formally unramified adic morphism of formal schemes locally of formal finite presentation.
    Then the infinitesimal neighbourhood
    \begin{align*}
        X \xrightarrow{\hat{f}} \hat{Y} \xrightarrow{\iota} Y
    \end{align*}
    of $f$ exists.
    Furthermore, $\hat{f}$ is a thickening locally of formal finite presentation and $\iota$ is a formally étale morphism locally of formal finite presentation.
\end{mainthm}

In particular, in the category of formal schemes, a formally unramified adic morphism locally of formal finite presentation can be expressed as the composition of a closed immersion and a formally étale morphism. \par

Note that an immersion of schemes locally of formal finite presentation is formally unramified and adic.
In this case, the infinitesimal neighbourhood is locally obtained by completing along the kernel ideal sheaf of $f : X \rightarrow Y$ (Example \ref{ex:immersions}). \par

Closely related results are obtained in \cite{MR2313672} and \cite{MR2497584} for immersions of locally Noetherian formal schemes. \par

\subsection*{Outline}
In \S \ref{sec:formal_schemes} we revise adic rings and formal schemes by recalling some basic results on tensor products and localisations.
In \S \ref{sec:closed_immersions} we define closed immersions of formal schemes as locally described by surjective adic morphisms.
In \S \ref{sec:thickenings} we define thickenings of formal schemes as closed immersions inducing a homeomorphism of topological spaces.
In \S \ref{sec:formal_finite_presentation} we define the property of being locally of formal finite presentation and show that it is a local property on both source and target.
In \S \ref{sec:infinitesimal_neighbourhoods} we construct the infinitesimal neighbourhood of formally unramified morphisms and show that the universal property holds.

\subsection*{Acknowledgements}
I would like to thank my PhD advisor Paolo Cascini for introducing me to foliations.
I also thank Richard Thomas and Calum Spicer for interest in the topic and for several corrections.
I thank Lambert A'Campo, Alessio Bottini, George Boxer, Riccardo Carini, Przemysław Grabowski and Johannes Nicaise for numerous useful conversations about formal schemes.
Finally, I thank the anonymous referee for several corrections and extensive comments and suggestions.
This research was conducted at the Max-Planck-Institut für Mathematik and Imperial College London, which I thank for their hospitality and financial support.

%% file: files/section1.tex
\section{Formal schemes}
\label{sec:formal_schemes}

\subsection{Adic rings}
\label{subsec:formal_schemes_local}

We begin with some revision and conventions. \par

In this paper, by an \emph{adic ring}, we mean a complete topological ring $A$ whose topology is the $I$-adic topology, for some \emph{finitely generated} ideal $I$ of $A$.
The ideal $I$ is said to be an \emph{ideal of definition}.
By convention, any ideal of definition is assumed to be finitely generated.

It is easy to check that an adic ring $A$ with ideal of definition $I$ is isomorphic to
\begin{align}
    \label{eq:filtered_limit}
    \lim{n \in \mathbb{N}} A / I^{n+1},
\end{align}
a filtered limit of discrete rings. \par

A \emph{morphism of adic rings} $\varphi : A \rightarrow B$ is a morphism of rings which is continuous as a morphism of topological spaces.
Any such morphism can be naturally viewed as an element of 
\begin{align}
    \label{eq:morphisms_adic_rings}
    \lim{m \in \mathbb{N}} \left( \colim{n \in \mathbb{N}} \left( \mathrm{Hom}_{\mathrm{Rings}} \left( A_n, B_m \right)\right) \right),
\end{align}
where $A_n = A / I^{n+1}$ and $B_m = B / J^{m+1}$ for some ideals of definition $I \subseteq A$ and $J \subseteq B$. \par

An \emph{adic morphism} $\varphi : A \rightarrow B$ is a morphism of adic rings such that for all ideals of definition $I \subseteq A$, the ideal extension $\varphi(I) \cdot A$ is an ideal of definition of $B$.
It suffices to check this condition for some ideal of definition of $A$ (see \cite[\href{https://stacks.math.columbia.edu/tag/0GBR}{Definition 0GBR}]{stacks-project} and the paragraph immediately following).
The composition of two adic morphisms is again an adic morphism (\cite[\href{https://stacks.math.columbia.edu/tag/0GXA}{Lemma 0GXA}]{stacks-project}). \par

Given morphisms of adic rings $\varphi : A \rightarrow B$ and $\psi : A \rightarrow A^{\prime}$, we may construct the \emph{completed tensor product} $B^{\prime} := B \hat{\otimes}_A A^{\prime}$ by completing the tensor product $B \otimes_A A^{\prime}$ with respect to the induced topology.
Since the ideals of definition are assumed to be finitely generated, $B^{\prime}$ is an adic ring (part (4) of \cite[\href{https://stacks.math.columbia.edu/tag/0GB4}{Lemma 0GB4}]{stacks-project}).
There are induced morphisms yielding a co-Cartesian square
\begin{equation}
    \label{diag:tensor_product_adic_rings}
	\begin{tikzcd}
		A \arrow[d, "\psi"] \arrow[r, "\varphi"] & B \arrow[d, "\psi^{\prime}"] \\
		A^{\prime} \arrow[r, "\varphi^{\prime}"] & B^{\prime}
	\end{tikzcd}
\end{equation}
in the category of adic rings.
Furthermore, if $I^{\prime}$ and $J$ are ideals of definition of $A^{\prime}$ and $B$ respectively, then $I^{\prime} \cdot B^{\prime} + J \cdot B^{\prime}$ is an ideal of definition of $B^{\prime}$. \par

If $\varphi : A \rightarrow B$ is an adic morphism and $\psi : A \rightarrow A^{\prime}$ is a morphism of adic rings, then the base change $\varphi^{\prime} : A^{\prime} \rightarrow B \hat{\otimes}_A A^{\prime}$ is an adic morphism (combine \cite[\href{https://stacks.math.columbia.edu/tag/0GX5}{Lemma 0GX5}]{stacks-project} and \cite[\href{https://stacks.math.columbia.edu/tag/0APU}{Lemma 0APU}]{stacks-project}). \par

If $A$ is an adic ring and $f \in A$ is an element, then the \emph{completed localisation} $A_{\{f\}}$ is the completion of the localisation $A_f$ with respect to the induced topology.
The morphism $A \rightarrow A_{\{f\}}$ is initial amongst all morphisms of adic rings $A \rightarrow B$ where $f$ invertible in $B$.

Given an adic ring $A$ and an element $f \in A$, the localisation morphism $A \rightarrow A_{\{f\}}$ is an adic morphism.
Furthermore, for any morphism $\varphi : A \rightarrow B$ of adic rings, there exists an induced morphism $\varphi_{\{f\}} : A_{\{f\}} \rightarrow B_{\{\varphi(f)\}}$ and a co-Cartesian diagram
\begin{equation}
    \label{diag:localisation_adic_rings}
	\begin{tikzcd}
		A \arrow[d] \arrow[r, "\varphi"] & B \arrow[d] \\
		A_{\{f\}} \arrow[r, "\varphi_{\{f\}}"] & B_{\{\varphi(f)\}}.
	\end{tikzcd}
\end{equation}
This fact follows from noting that the universal property of the completed localisation applied to $B_{\{\varphi(f)\}}$ is the universal property of the completed tensor product $B \hat{\otimes}_A A_{\{f\}}$.

Being adic is a local property.
More precisely,
\begin{enumerate}
	\item if $A$ is an adic ring with ideal of definition $I$, then for any $f \in A$, $A_{\{f\}}$ is an adic ring with ideal of definition $I_{\{f\}}$ (part (ii) of \cite[Proposition (7.6.11)]{MR217083});
	\item if $A$ is a complete topological ring endowed with a linear topology and there exist elements $f_1, \ldots, f_r$ generating the unit ideal of $A$ such that $A_{\{f_1\}}, \ldots, A_{\{f_r\}}$ are adic rings, then $A$ is an adic ring (special case of \cite[\href{https://stacks.math.columbia.edu/tag/0GXQ}{Lemma 0GXQ}]{stacks-project}).
\end{enumerate}

Given an adic ring $A$ with ideal of definition $I$, define the \emph{reduction}
\begin{align}
    \label{eq:reduction_affine}
    A_{\mathrm{red}} = A / \left( \mathrm{rad}_A \, I \right).
\end{align}
Since $I$ is an ideal of definition, $\mathrm{rad}_A \, I$ is equal to the \emph{topological nilradical} of $A$, the largest ideal of definition of $A$.
It follows that $A_{\mathrm{red}}$ is a reduced discrete ring and the morphism $A \rightarrow A_{\mathrm{red}}$ is initial amongst all morphisms of adic rings $A \rightarrow B$ where $B$ is a reduced discrete ring.
Given a morphism of adic rings $\varphi : A \rightarrow B$, there exists a unique induced morphism of discrete rings $\varphi_{\mathrm{red}} : A_{\mathrm{red}} \rightarrow B_{\mathrm{red}}$.
Thus, reduction of a ring is the left adjoint to the inclusion of discrete reduced rings into adic rings, hence it commutes with colimits. \par

\subsection{Formal schemes}
\label{subsec:formal_schemes_global}

In this paper, by a \emph{formal scheme}, we mean a locally topologically ringed space $X$ (see \cite[\href{https://stacks.math.columbia.edu/tag/0AHY}{Section 0AHY}]{stacks-project}) which is locally isomorphic to the formal spectrum of an adic ring (with a finitely generated ideal of definition). \par

Given a formal scheme $X$, the affine formal open subschemes of $X$ form a basis for its topology.
Indeed, by \cite[Proposition (10.4.3)]{MR217083}, there is a basis of open subschemes of $X$ which are formal spectra of complete topological rings.
But since being adic is a local property (\S \ref{subsec:formal_schemes_local}), such complete topological rings must in fact be adic, thus such basis consists of affine formal open subschemes.

A formal scheme $X$ is locally Noetherian if every affine open subscheme is the formal spectrum of an adic Noetherian ring.
This can be checked on an affine open covering (\cite[\href{https://stacks.math.columbia.edu/tag/0AKX}{Lemma 0AKX}]{stacks-project}). \par

A morphism of formal schemes is a morphism of locally topologically ringed spaces.
A morphism of formal schemes is \emph{adic} if for any affine open set $V = \spf A \subseteq Y$ and $U = \spf B \subseteq f^{-1}(V) \subseteq X$, the induced morphism of adic rings $A \rightarrow B$ is adic.
The composition of adic morphisms is again adic. \par

The fibre product of formal schemes exist.
More precisely, if $f : X \rightarrow S$ and $g : Y \rightarrow S$ are morphisms of formal schemes, the fibre product $X \times_S Y$ is a formal scheme.
The existence can be shown as in \cite[Proposition (10.7.2)]{MR217083} on using the fact that the tensor product of adic rings is adic. \par

Given a formal scheme $X$, there exists a reduced scheme $X_{\mathrm{red}} \rightarrow X$ which is terminal amongst all morphisms of formal schemes $Z \rightarrow X$, where $Z$ is a reduced scheme (see also \cite[\href{https://stacks.math.columbia.edu/tag/0AIN}{Lemma 0AIN}]{stacks-project}).
This is the \emph{reduction} of $X$.
We can construct $X_{\mathrm{red}} \rightarrow X$ locally on spectra using (\ref{eq:reduction_affine}), and then glue using the fact that completed localisation is a colimit and reduction of a ring commutes with colimits.

%% file: files/section2.tex
\section{Closed immersions}
\label{sec:closed_immersions}

\subsection{Surjective adic morphisms}
\label{subsec:closed_immersions_local}

The affine local model of a closed immersion of formal schemes is a surjection of adic rings $\varphi : A \rightarrow B$ (with finitely generated ideals of definition) such that $\varphi$ is an adic morphism.
The latter condition is equivalent to $\ker \varphi$ being a closed ideal of $A$, and is also equivalent to $\varphi$ inducing an isomorphism of adic rings $A / \left( \ker \varphi \right) = B$.

\begin{example}
    \label{ex:diagonal_closed_immersion}
	Let $A \rightarrow B$ be a morphism of adic rings.
	Then the diagonal morphism
	\begin{align}
	    \label{eq:diagonal_morphism}
		B \hat{\otimes}_A B &\rightarrow B \\
		b \otimes b^{\prime} &\rightarrow b b^{\prime} \nonumber
	\end{align}
	is a surjective adic morphism of adic rings.
	To see this, suppose that $J \subseteq B$ is an ideal of definition of $B$ and recall that $B \hat{\otimes}_A B$ is an adic ring with ideal of definition given by the image of 
    \begin{align}
	    \label{eq:ideal_definition_diagonal}
		J \otimes_A B + B \otimes_A J.
	\end{align}
	Certainly, (\ref{eq:diagonal_morphism}) is surjective.
	Furthermore, the image of the ideal of definition in (\ref{eq:ideal_definition_diagonal}) under the diagonal morphism is $J \subseteq B$, an ideal of definition of $B$.
	Therefore (\ref{eq:diagonal_morphism}) is an adic morphism.
\end{example}

We firstly show that being a surjective adic morphism is preserved by base change.

\begin{lemma}
    \label{lem:base_change_surjective_adic}
	The base change of a surjective adic morphism is a surjective adic morphism.
\end{lemma}

\begin{proof}
    Suppose $\varphi : A \rightarrow B$ is a surjective adic morphism and let $\psi: A \rightarrow A^{\prime}$ be a morphism of adic rings.
    Let
    \begin{align}
        \label{eq:base_change_surjective_adic}
        \varphi^{\prime} : A^{\prime} \rightarrow A^{\prime} \hat{\otimes}_A B =: B^{\prime}
    \end{align}
    denote the base change of $\varphi$.
    Then $B^{\prime}$ is an adic ring and $\varphi^{\prime}$, being the base change of $\varphi$, is an adic morphism.
    It remains to show that $\varphi^{\prime}$ is surjective.
	Let $I^{\prime}$ be an ideal of definition of $A^{\prime}$.
	By part (1) of \cite[\href{https://stacks.math.columbia.edu/tag/0315}{Lemma 0315}]{stacks-project}, it suffices to show that the induced morphism $A^{\prime} / I^{\prime} \rightarrow B^{\prime} / (I^{\prime} \cdot B^{\prime})$ is surjective.
	But the latter morphism is the base change of $\varphi$ by the morphism $A \rightarrow A^{\prime} / I^{\prime}$.
	Since $\varphi$ is surjective, the lemma follows.
\end{proof}

It is clear that the composition of surjective adic morphisms is a surjective adic morphism. \par

Next, we show that being a surjective adic morphism can be checked on localisations of the domain.

\begin{lemma}
    \label{lem:surjective_adic_local}
	Let $\varphi : A \rightarrow B$ be a morphism of adic rings.
	\begin{enumerate}
		\item If $\varphi$ is a surjective adic morphism, then, for any $f \in A$, the complete localisation morphism
		\begin{align}
		    \label{eq:surjective_adic_localisation_i}
			\varphi_{\{f\}} : A_{\{f\}} \rightarrow B_{\{\varphi(f)\}}
		\end{align}
		is a surjective adic morphism, and,
		\item if there exist elements $f_1, \ldots, f_r$ in $A$ generating the unit ideal such that the complete localisation morphism
		\begin{align}
		    \label{eq:surjective_adic_localisation_ii}
			\varphi_{\{f_i\}} : A_{\{f_i\}} \rightarrow B_{\{\varphi(f_i)\}}
		\end{align}
		is a surjective adic morphism for all $i \leq r$, $\varphi$ is a surjective adic morphism.
	\end{enumerate}
\end{lemma}

\begin{proof}
	For (1), since $\varphi_{\{f\}}$ is the base change of $\varphi$, the claim follows from Lemma \ref{lem:base_change_surjective_adic}. \par
	
	For (2), we first show that $\varphi$ is adic.
	Let $I$ and $J$ be ideals of definition of $A$ and $B$ respectively.
	Furthermore, let $J^{\prime} = I \cdot B$.
	Fix $i \leq r$.
	Since $I_{\{f_i\}}$ and $J_{\{f_i\}}^{\prime}$ are the extensions of the ideals $I$ and $J^{\prime}$ respectively (\cite[(7.6.9)]{MR217083}), it follows that $J_{\{f_i\}}^{\prime} = I_{\{f_i\}} \cdot B_{\{f_i\}}$.
	But now $I_{\{f_i\}}$ is an ideal of definition of $A_{\{f_i\}}$ and since, by assumption, $\varphi_{\{f_i\}}$ is adic, $J_{\{f_i\}}^{\prime}$ is an ideal of definition of $B_{\{f_i\}}$.
	Similarly, $J_{\{f_i\}}$ is an ideal of definition of $B_{\{f_i\}}$, therefore $J_{\{f_i\}}$ and $J_{\{f_i\}}^{\prime}$ define the same adic topology on $B_{\{f_i\}}$.
	Since containment of ideals is a local property, a simple computation shows that $J$ and $J^{\prime}$ define the same adic topology on $B$. \par 
	
	Next we show that $\varphi$ is surjective.
	Since $\varphi$ is an adic morphism, part (1) of \cite[\href{https://stacks.math.columbia.edu/tag/0315}{Lemma 0315}]{stacks-project} allows us to reduce to the case when both $A$ and $B$ are discrete rings.
	In this case, it is well-known that being surjective is an affine-local property.
\end{proof}

\subsection{Closed immersions of formal schemes}
\label{subsec:closed_immersions_global}

\begin{definition}
    \label{def:closed_immersion}
	A morphism of formal schemes $f : X \rightarrow Y$ is a \emph{closed immersion} if for any affine open formal subscheme $V = \spf A \subseteq Y$, the pre-image $U := f^{-1}(V) = \spf B$ is affine and the induced morphism $A \rightarrow B$ is a surjective adic morphism.
    An \emph{immersion} $f : X \rightarrow Y$ of formal schemes is a closed immersion into an open subscheme of $Y$.
\end{definition}

Note that immersions are adic morphisms.
Indeed, closed immersions are adic by definition and open immersions are adic since so is the morphism $A \rightarrow A_{\{f\}}$ for any $f \in A$. \par

We firstly show that being a closed immersion is an affine-local property.

\begin{lemma}
    \label{lem:closed_immersion_local}
	Suppose that $f : X \rightarrow Y$ is a morphism of formal schemes and assume there exists an affine open covering $\{V_i = \spf A_i\}_{i \in \Lambda}$ of $Y$ such that the pre-image $U_i := f^{-1}(V_i) = \spf B_i$ is affine and the induced morphism $A_i \rightarrow B_i$ is a surjective adic morphism.
	Then $f$ is a closed immersion of formal schemes.
\end{lemma}

\begin{proof}
	Let $V = \spf A \subseteq Y$ be an affine open formal  subscheme.
	It suffices to show that $U := f^{-1}(V) \subseteq X$ is an affine open subscheme.
	Indeed, we can reduce to Lemma \ref{lem:surjective_adic_local} by employing the fact that $Y$ has a basis of affine open formal subschemes. \par

	Without loss of generality, we may assume $V = Y$ and prove that $X$ is affine.
	To this end, it suffices to prove its reduction $X_{\mathrm{red}}$ is affine (necessary implication of \cite[\href{https://stacks.math.columbia.edu/tag/0DE8}{Lemma 0DE8}]{stacks-project}, based on \cite[Proposition 3.32]{MR2520785}).
	By the sufficient implication of the same result, $Y_{\mathrm{red}}$ is affine.
	But now the assumption on $f$ implies that $f$ is an affine morphism of schemes (\cite[\href{https://stacks.math.columbia.edu/tag/01S8}{Lemma 01S8}]{stacks-project}), and we deduce that $X_{\mathrm{red}}$ is affine. 
\end{proof}

In particular, $\varphi : A \rightarrow B$ is a surjective adic morphism if and only if the induced morphism $\spf \varphi : \spf B \rightarrow \spf A$ is a closed immersion of formal schemes.

\begin{lemma}
    \label{lem:base_change_immersions}
	The base change of a closed immersion of formal schemes is a closed immersion of formal schemes.
    The base change of an immersion of formal schemes is an immersion of formal schemes.
\end{lemma}

\begin{proof}
	Let $f : X \rightarrow Y$ be a closed immersion and let $g : Y^{\prime} \rightarrow Y$ a morphism of formal schemes.
	Let $f^{\prime} : X^{\prime} = X \times_Y Y^{\prime} \rightarrow Y$ denote the base change.
	By Lemma \ref{lem:closed_immersion_local}, we may check that $f^{\prime}$ is a closed immersion on an affine open covering of $Y^{\prime}$, hence $Y$ and $Y^{\prime}$ can be assumed affine.
    Since $f$ is a closed immersion, $X$ can also be assumed affine, hence the result follows from Lemma \ref{lem:base_change_surjective_adic}. \par

    The base change of an open immersion of formal schemes is certainly an open immersion of formal schemes, hence the same is true for an immersion of formal schemes.
\end{proof}

Next, we check that the diagonal of a morphism is an immersion.
This allows us to discuss separation properties of morphisms.

\begin{lemma}
    \label{lem:diagonal_immersion}
    Let $f : X \rightarrow Y$ be a morphism of formal schemes.
    Then the diagonal morphism $\Delta_f : X \rightarrow X \times_Y X$ is an immersion of formal schemes.
\end{lemma}

\begin{proof}
    As in \cite[\href{https://stacks.math.columbia.edu/tag/01KJ}{Lemma 01KJ}]{stacks-project}, the claim readily reduces to Example \ref{ex:diagonal_closed_immersion}.
\end{proof}

It is clear that the composition of closed immersions and immersions of formal schemes is a closed immersion and an immersion of formal schemes respectively.
Indeed, for the latter case, note that the same proof as in \cite[\href{https://stacks.math.columbia.edu/tag/02V0}{Lemma 02V0}]{stacks-project} works in the case of formal schemes.
Furthermore, immersions of formal schemes are monomorphisms, hence separated.

\begin{lemma}
    \label{lem:closed_immersion_factor}
    Let 
    \begin{align}
        X \xrightarrow{f} Y \xrightarrow{g} Z
	\end{align}
	be morphisms of formal schemes.
	If $g \circ f$ is an immersion, so is $f$.
\end{lemma}

\begin{proof}
    For schemes, this is part (1) of \cite[\href{https://stacks.math.columbia.edu/tag/07RK}{Lemma 07RK}]{stacks-project}.
    Note however that the proof also works for formal schemes.
    Indeed, the proof relies on the fact that the composition of two immersions is an immersion, the fact that the base change of an immersion is an immersion (Lemma \ref{lem:base_change_immersions}) and on \cite[\href{https://stacks.math.columbia.edu/tag/01KR}{Lemma 01KR}]{stacks-project}.
    The latter also holds for formal schemes, since it uses general category theory and the fact that the diagonal morphism is an immersion (Lemma \ref{lem:diagonal_immersion}).
\end{proof}

%% file: files/section3.tex
\section{Thickenings}
\label{sec:thickenings}

\subsection{Thickenings of adic rings}
\label{subsec:thickenings_local}

Analogously to the case of schemes, a thickening of formal schemes is defined as a closed immersion which induces a homeomorphism of topological spaces.

\begin{definition}
A morphism of adic rings $\varphi : A \rightarrow B$ is a \emph{thickening} if it is a surjective adic morphism such that
\begin{align}
    \label{eq:thickening_formal_spectra}
	\spf \varphi : \spf B \rightarrow \spf A
\end{align}
induces a homeomorphism of topological spaces.
\end{definition}

Recall that a surjective morphism of (discrete) rings is a thickening if and only if every element of its kernel is nilpotent.
A similar characterisation holds for adic rings.

\begin{lemma}
    \label{lem:characterisation_thickenings}
	Let $\varphi : A \rightarrow B$ be a surjective adic morphism of adic rings and let $I$ and $J$ be ideals of definition of $A$ and $B$ respectively.
	Then the following are equivalent:
	\begin{enumerate}
		\item $\varphi$ is a thickening;
		\item $\ker \varphi \subseteq \mathrm{rad}_{A} \, I$;
		\item $\mathrm{rad}_{A} \, I = \varphi^{-1} \left( \mathrm{rad}_{B} \, J \right)$.
	\end{enumerate}
\end{lemma}

\begin{proof}
	$(1) \leftrightarrow (2)$.
	Since $\varphi$ is adic, $I \cdot B$ is an ideal of definition of $B$, hence there is an induced surjective morphism of discrete rings
	\begin{align}
		\varphi_I : A / I \rightarrow B / \left( I \cdot B \right).
	\end{align}
	By construction, $\varphi$ is a thickening of adic rings if and only if $\varphi_I$ is a thickening of rings.
	By standard algebra, $\varphi_I$ is a thickening of rings if and only if $\ker \left( \varphi_I \right) \subseteq \mathrm{rad}_{A/I} \, 0$.
	This is equivalent to requiring $\ker \varphi \subseteq \mathrm{rad}_{A} \, I$.
	
	$(2) \leftrightarrow (3)$.
	Recall that $\mathrm{rad}_{B} \, J$ is equal to the topological nilradical of the adic ring $B$ and hence is independent of the ideal of definition chosen.
	We may therefore take $J = I \cdot B$.
	Since $\varphi$ is continuous, $(3)$ is equivalent to
	\begin{align}
	    \label{eq:reverse_inclusion_thickening}
		\varphi^{-1} \left( \mathrm{rad}_{B} \, J \right) \subseteq \mathrm{rad}_{A} \, I.
	\end{align}
	By construction, $\varphi^{-1} \left( J \right) = I + \ker \varphi$, hence (\ref{eq:reverse_inclusion_thickening}) holds if and only if $\ker \varphi \subseteq \mathrm{rad}_A \, I$.
\end{proof}

\begin{lemma}
    \label{lem:base_change_thickening_affine}
	The base change of a thickening of adic rings is a thickening of adic rings.
\end{lemma}

\begin{proof}
    Suppose $\varphi : A \rightarrow B$ is a thickening of adic rings and let $\psi: A \rightarrow A^{\prime}$ be a morphism of adic rings.
    Let
    \begin{align}
        \label{eq:base_change_thickening}
        \varphi^{\prime} : A^{\prime} \rightarrow A^{\prime} \hat{\otimes}_A B =: B^{\prime}
    \end{align}
    denote the base change of $\varphi$.
	By Lemma \ref{lem:base_change_surjective_adic}, $\varphi^{\prime}$ is a surjective adic morphism, hence it suffices to show it induces a homeomorphism of topological spaces.
	Let $I$ and $I^{\prime}$ be ideals of definition of $A$ and $A^{\prime}$ respectively and let $K$ and $K^{\prime}$ be the kernel ideals of $\varphi$ and $\varphi^{\prime}$ respectively.
	In order to show $\varphi^{\prime}$ is a thickening, by the $(2) \rightarrow (1)$ implication of Lemma \ref{lem:characterisation_thickenings}, it suffices to show that $K^{\prime} \subseteq \mathrm{rad}_{A^{\prime}} \, I^{\prime}$.
	Up to replacing $I$ by a sufficiently high enough power, we can assume that
	\begin{align}
	    \label{eq:psi_continuous}
		I \subseteq \psi^{-1} \left( I^{\prime} \right).
	\end{align}
	Let $x^{\prime} \in K^{\prime}$.
	By assumption, $K^{\prime} = K \cdot A^{\prime}$, hence $x^{\prime}$ can be written as $\sum_{\lambda \in \Lambda} a_{\lambda}^{\prime} x_{\lambda}$, where $a_{\lambda}^{\prime} \in A^{\prime}$, $x_{\lambda} \in K$ and $\Lambda$ is a finite set.
	Since $\varphi$ is a thickening, by the $(1) \rightarrow (2)$ implication of Lemma \ref{lem:characterisation_thickenings}, for all $\lambda \in \Lambda$, there exists an $n_{\lambda}$ such that $x_{\lambda}^{n_{\lambda}+1} \in I$.
	Set $n := \sum_{\lambda \in \Lambda} n_{\lambda}$, then it is straightforward to see that $\left( x^{\prime} \right)^{n+1} \in I \cdot A^{\prime}$.
	Now (\ref{eq:psi_continuous}) yields that $x^{\prime} \in \mathrm{rad}_{A^{\prime}} \, I^{\prime}$.
\end{proof}

It is clear that the composition of thickenings is a thickening.

\begin{lemma}
    \label{lem:thickening_affine_local}
	Part (1) and part (2) of Lemma \ref{lem:surjective_adic_local} for surjective adic morphisms are also true for thickenings of adic rings.
\end{lemma}

\begin{proof}
	For (1), since $\varphi_{\{f\}}$ is the base change of $\varphi$, the claim follows from Lemma \ref{lem:base_change_thickening_affine}. \par
	
	For (2), by Lemma \ref{lem:surjective_adic_local}, it suffices to show that $\varphi$ induces a homeomorphism of topological spaces.
	Let $I$ be an ideal of definition of $A$.
	By the $(2) \rightarrow (1)$ implication of Lemma \ref{lem:characterisation_thickenings}, it suffices to show that $\ker \varphi \subseteq \mathrm{rad}_A \, I$.
	This is true since containment of ideals is a local property.
\end{proof}

\subsection{Thickenings of formal schemes}
\label{subsec:thickenings_global}

\begin{definition}
    \label{def:thickening}
	A morphism of formal schemes $f : X \rightarrow Y$ is a \emph{thickening} if for any affine open formal subscheme $V = \spf A \subseteq Y$, the pre-image $U := f^{-1}(V) = \spf B$ is affine and the induced morphism $A \rightarrow B$ is a thickening of adic rings.
\end{definition}

Since being a homeomorphism is local on the target, it is easy to see that a morphism of formal schemes is a thickening if and only if it is a closed immersion of formal schemes inducing a homeomorphism of underlying topological spaces.
It follows that $\varphi : A \rightarrow B$ is a thickening of adic rings if and only if the induced morphism $\spf \varphi : \spf B \rightarrow \spf A$ is a thickening of formal schemes. \par

Being a thickening is an affine-local property.

\begin{lemma}
    \label{lem:thickening_local}
	Suppose that $f : X \rightarrow Y$ is a morphism of formal schemes and assume there exists an affine open covering $\{V_i = \spf A_i\}_{i \in \Lambda}$ of $Y$ such that the pre-image $U_i := f^{-1}(V_i) = \spf B_i$ is affine and the induced morphism $A_i \rightarrow B_i$ is a thickening of adic rings.
	Then $f$ is a thickening of formal schemes.
\end{lemma}

\begin{proof}
	As in the proof of Lemma \ref{lem:closed_immersion_local}, we can reduce to the corresponding local statement (Lemma \ref{lem:thickening_affine_local}).
\end{proof}

\begin{lemma}
    \label{lem:base_change_thickenings}
	The base change of a thickening of formal schemes is a thickening of formal schemes.
\end{lemma}

\begin{proof}
	As in the proof of Lemma \ref{lem:base_change_immersions}, we can reduce to the corresponding local statement (Lemma \ref{lem:base_change_thickening_affine}).
\end{proof}

It is clear that the composition of thickenings of formal schemes is a thickening of formal schemes.

\begin{lemma}
    \label{lem:image_thickening}
    Let $f : X \rightarrow Y$ be a thickening of formal schemes.
    If $X$ is affine, so is $Y$.
\end{lemma}

\begin{proof}
    Note that $f$ induces an isomorphism $f_{\mathrm{red}} : X_{\mathrm{red}} \rightarrow Y_{\mathrm{red}}$.
	Indeed, this can be checked affine locally on $Y$, where the claim follows from the $(1) \rightarrow (3)$ implication of Lemma \ref{lem:characterisation_thickenings} and the definition of reduction.
	We conclude by noting that a formal scheme is affine if and only if so is its reduction (\cite[\href{https://stacks.math.columbia.edu/tag/0DE8}{Lemma 0DE8}]{stacks-project}).
\end{proof}

%% file: files/section4.tex
\section{Formal finite presentation}
\label{sec:formal_finite_presentation}

\subsection{Morphisms of adic rings}
\label{subsec:presentation_local}

We give a definition of finite presentation which takes the topologies of the rings into account.

\begin{lemma}
    \label{lem:finite_presentation}
    Let $\varphi : A \rightarrow B$ be a morphism of adic rings, then the following conditions are equivalent.
    \begin{enumerate}
        \item[(1A)] For all ideals of definition $I \subseteq A$ and $J \subseteq B$ and for all $m \in \mathbb{N}$, there exists an $n \in \mathbb{N}$ such that the induced morphism $\varphi_{n,m} : A / I^{n+1} \rightarrow B / J^{m+1}$ is of finite presentation.
        \item[(1B)] There exist ideals of definition $I \subseteq A$ and $J \subseteq B$ with the following property: for all $m \in \mathbb{N}$, there exists an $n \in \mathbb{N}$ such that the induced morphism $\varphi_{n,m} : A / I^{n+1} \rightarrow B / J^{m+1}$ is of finite presentation.
        \item[(2A)] For all ideals of definition $J \subseteq B$ and for all $m \in \mathbb{N}$, the induced morphism $\varphi_{m} : A \rightarrow B / J^{m+1}$ is of finite presentation.
        \item[(2B)] There exists an ideal of definition $J \subseteq B$ such that for all $m \in \mathbb{N}$, the induced morphism $\varphi_{m} : A \rightarrow B / J^{m+1}$ is of finite presentation.
    \end{enumerate}
\end{lemma}

\begin{proof}
    It is clear that $(1A) \rightarrow (1B)$ and $(2A) \rightarrow (2B)$.
    The implications $(1A) \rightarrow (2A)$ and $(1B) \rightarrow (2B)$ follow from the fact that $A \rightarrow A / I^{n+1}$ is of finite presentation.
    The implications $(2A) \rightarrow (1A)$ and $(2B) \rightarrow (1B)$ follow from part (4) of \cite[\href{https://stacks.math.columbia.edu/tag/00F4}{Lemma 00F4}]{stacks-project}. \par
    
    It remains to show that $(2B) \rightarrow (2A)$.
    Let $J$ and $J^{\prime}$ be two ideals of definition of $B$ and assume $J^{\prime}$ satisfies $(2B)$.
    For any $m \in \mathbb{N}$, $J^m$ is open, hence there exists an $m^{\prime} \in \mathbb{N}$ such that $\left(J^{\prime} \right)^{m^{\prime}} \subseteq J^m$.
    Since $J^m$ is finitely generated, the induced morphism $B / \left(J^{\prime} \right)^{m^{\prime}} \rightarrow B / J^m$ is of finite presentation and the result follows.
\end{proof}

\begin{definition}
    \label{def:finite_presentation_affine}
    A morphism $\varphi : A \rightarrow B$ of adic rings is of \emph{formal finite presentation} if any of the equivalent conditions of Lemma \ref{lem:finite_presentation} are satisfied.
\end{definition}

\begin{lemma}
    \label{lem:base_change_finite_presentation_affine}
	The base change of a morphism of formal finite presentation is a morphism of formal finite presentation.
\end{lemma}

\begin{proof}
	Using the characterisations in Lemma \ref{lem:finite_presentation}, it is straightforward to reduce to the case of discrete rings.
	Then the claim is true by part (4) of \cite[\href{https://stacks.math.columbia.edu/tag/05G5}{Lemma 05G5}]{stacks-project}.
\end{proof}

It is clear that the composition of morphisms of formal finite presentation is a morphism of formal finite presentation.

\begin{lemma}
    \label{lem:localisation_finite_presentation}
    Let $A$ be an adic ring and let $f \in A$.
    Then the completed localisation morphism $A \rightarrow A_{\{f\}}$ is of formal finite presentation.
\end{lemma}

\begin{proof}
    This is true since the localisation $A_f$ is of finite presentation over $A$ and the ideal of definition $I_{\{f\}} \subseteq A_{\{f\}}$ is finitely generated.
\end{proof}

\begin{lemma}
    \label{lem:finite_presentation_affine_local}
	Being a morphism of formal finite presentation is a local property (see part (1) of \cite[\href{https://stacks.math.columbia.edu/tag/01SR}{Definition 01SR}]{stacks-project}).
\end{lemma}

\begin{proof}
    We need to show that parts (a), (b) and (c) of the definition are satisfied. \par

    Part (a) states that, if $A \rightarrow B$ is a morphism of formal finite presentation, so is $A_{\{f\}} \rightarrow B_{\{f\}}$ for any $f \in A$.
    This is true by Lemma \ref{lem:base_change_finite_presentation_affine}. \par

    Part (b) states that, if $A$ and $B$ are (adic) rings, $f \in A$ and $g \in B$ are elements, and $A_{\{f\}} \rightarrow B$ is of formal finite presentation, then so is $A \rightarrow B_{\{g\}}$.
    This is true by applying twice Lemma \ref{lem:localisation_finite_presentation} to the composition $A \rightarrow A_{\{f\}} \rightarrow B \rightarrow B_{\{g\}}$. \par
	
	Part (c) states that, if $A \rightarrow B$ is a morphism of adic rings and $A \rightarrow B_{\{g_i\}}$ is of formal finite presentation for $g_1, \ldots, g_r$ in $B$ generating the unit ideal, then $A \rightarrow B$ is of formal finite presentation.
    This easily follows from reduction to the discrete case (\cite[\href{https://stacks.math.columbia.edu/tag/00EP}{Lemma 00EP}]{stacks-project}).
\end{proof}

\begin{lemma}
    \label{lem:finite_presentation_factor_affine}
    Let
    \begin{align}
        \label{eq:composition_finite_presentation_affine}
        A \xrightarrow{\varphi} B \xrightarrow{\psi} C
    \end{align}
    be morphisms of adic rings.
    If $\varphi$ and $\psi \circ \varphi$ are of formal finite presentation, so is $\psi$.
\end{lemma}

\begin{proof}
    It is straightforward to reduce to the discrete case where the claim follows from part (4) of \cite[\href{https://stacks.math.columbia.edu/tag/00F4}{Lemma 00F4}]{stacks-project}.
\end{proof}

\begin{lemma}
    \label{lem:noetherian_ascent_affine}
    Let $\varphi : A \rightarrow B$ be a morphism of formal finite presentation.
    If $A$ is Noetherian, so is $B$.
\end{lemma}

\begin{proof}
    Let $J$ be an ideal of definition of $B$.
    By assumption, $J$ is finitely generated and $A \rightarrow B/J$ is of finite presentation.
    Since $A$ is Noetherian, so is $B/J$ (\cite[\href{https://stacks.math.columbia.edu/tag/00FN}{Lemma 00FN}]{stacks-project}).
    Now completeness of $B$ implies that $B$ is Noetherian (\cite[\href{https://stacks.math.columbia.edu/tag/05GH}{Lemma 05GH}]{stacks-project}).
\end{proof}

Next, we work towards showing that any formally smooth morphism which is also a thickening must be the identity morphism.

\begin{lemma}
    \label{lem:chinese_remainder}
    Let $R$ be a discrete ring and let $I$ and $J$ be ideals contained in the nilradical of $R$.
    Then the induced morphism
    \begin{align}
        \label{eq:chinese_remainder}
        R \rightarrow \left( R / I \right) \times_{R / (I+J)} \left( R / J \right)
    \end{align}
    is surjective and its kernel is contained in the nilradical of $R$.
\end{lemma}

\begin{proof}
    By the generalised Chinese Remainder Theorem, we have that $\left( R / I \right) \times_{R / (I+J)} \left( R / J \right) = R/(I \cap J)$, hence (\ref{eq:chinese_remainder}) is surjective.
    Furthermore, its kernel is $I \cap J$.
    Since $I$ and $J$ are contained in the nilradical of $R$, so is $I \cap J$. 
\end{proof}

\begin{lemma}
    \label{lem:formally_smooth_lift_affine}
    Let $\varphi : A \rightarrow B$ be a formally smooth morphism of formal finite presentation and let $\varphi^{\prime} : C^{\prime} \rightarrow C$ be a thickening of adic rings fitting in a commutative solid diagram
    \begin{equation}
        \label{diag:formally_smooth_lift}
        \begin{tikzcd}
            A \arrow[d, "\varphi"] \arrow[r, "\psi"] & C^{\prime} \arrow[d, "\varphi^{\prime}"] \\
            B \arrow[r, "\psi^{\prime}"] \arrow[ur, dashed] & C.
        \end{tikzcd}
    \end{equation}
    Then there exists a dashed morphism making the diagram commute.
\end{lemma}

\begin{proof}
    We firstly assume $\varphi^{\prime} : C^{\prime} \rightarrow C$ is a thickening of discrete rings.
    Let $J$ be an ideal of definition of $B$.
    Since $\psi^{\prime}$ is continuous and $C$ is discrete, $\psi^{\prime}$ factorises as $B \rightarrow B/J^{n+1} \rightarrow C$, for some $n \in \mathbb{N}$ large enough.
    By assumption, the composition $A \rightarrow B/J^{n+1}$ is of finite presentation.
    Let $K$ be the kernel of $\varphi^{\prime}$.
    We may express $K$ as the colimit of its finitely generated submodules $K_{\lambda}$, for $\lambda \in \Lambda$, where $\Lambda$ is a filtered set (see part (1) of \cite[\href{https://stacks.math.columbia.edu/tag/02CP}{Exercise 02CP}]{stacks-project}).\
    Since colimits commute with taking quotients, it follows that
    \begin{align}
        \label{eq:colimit_quotient}
        C = \frac{C^{\prime}}{K} = \frac{C^{\prime}}{\colim{\lambda \in \Lambda} K_{\lambda}} = \colim{\lambda \in \Lambda} \left( \frac{C^{\prime}}{K_{\lambda}} \right).
    \end{align}
    But now, since $A \rightarrow B/J^{n+1}$ is of finite presentation, using the characterisation in \cite[\href{https://stacks.math.columbia.edu/tag/00QO}{Lemma 00QO}]{stacks-project}, we see that
    \begin{align}
        \label{eq:finite_presentation_factorisation}
        \colim{\lambda \in \Lambda} \mathrm{Hom}_A \left( \frac{B}{J^{n+1}}, \frac{C^{\prime}}{K_{\lambda}} \right) \xrightarrow{\sim} \mathrm{Hom}_A \left( \frac{B}{J^{n+1}}, \colim{\lambda \in \Lambda} \frac{C^{\prime}}{K_{\lambda}} \right) = \mathrm{Hom}_A \left( \frac{B}{J^{n+1}}, C \right),
    \end{align}
    hence the dotted morphism in the solid commutative diagram
    \begin{equation}
        \label{diag:formally_smooth_lift_helper}
        \begin{tikzcd}
            A \arrow[dd, "\varphi"] \arrow[rr, "\psi"] & & C^{\prime} \arrow[d] \\
            & & C^{\prime}/K_{\lambda} \arrow[d] \\
            B \arrow[r] \arrow[uurr, dashed] & B/J^{n+1} \arrow[r] \arrow[ur, dotted] & C
        \end{tikzcd}
    \end{equation}
    exists for some $\lambda \in \Lambda$ large enough.
    But now $K_{\lambda}$ is finitely generated and is contained in the nilradical of $C^{\prime}$, since so is $K$ (Lemma \ref{lem:characterisation_thickenings} for discrete rings).
    It follows that $K_{\lambda}^{k+1} = 0$ for $k \gg 0$, and we may express the thickening $C^{\prime} \rightarrow C^{\prime}/K_{\lambda}$ as a composition of first order thickenings of discrete rings.
    Now the dashed morphism in Diagram (\ref{diag:formally_smooth_lift_helper}) can be constructed inductively using the defining property of formal smoothness for $\varphi$. \par

    Now assume $\varphi^{\prime} : C^{\prime} \rightarrow C$ is a thickening of adic rings and let $L^{\prime}$ be an ideal of definition of $C^{\prime}$.
    We construct a lift $B \rightarrow C^{\prime}$ by inductively finding compatible lifts $B \rightarrow C^{\prime} / (L^{\prime})^{n+1}$ in the commutative diagram
    \begin{equation}
        \label{diag:formally_smooth_lift_induction}
        \begin{tikzcd}
            A \arrow[d, "\varphi"] \arrow[r, "\psi"] & C^{\prime} \arrow[d, "\varphi^{\prime}"] \arrow[r] & \ldots \arrow[d] \arrow[r] & C^{\prime} / (L^{\prime})^{n+1} \arrow[d] \arrow[r] & C^{\prime} / (L^{\prime})^{n} \arrow[d] \arrow[r] & \ldots \arrow[d] \arrow[r] & C^{\prime} / L^{\prime} \arrow[d] \\
            B \arrow[r, "\psi^{\prime}"] & C \arrow[r] & \ldots \arrow[r] & C/L^{n+1} \arrow[r] & C/L^{n} \arrow[r] & \ldots \arrow[r] & C/L.
        \end{tikzcd}
    \end{equation}
    Since $C^{\prime}$ is defined as a limit of discrete rings, this is sufficient to show that a compatible dashed morphism in Diagram (\ref{diag:formally_smooth_lift}) exists. \par

    Since $\varphi^{\prime}$ is adic, $L := L^{\prime} \cdot C$ is an ideal of definition of $C$, so the base change $C^{\prime} / (L^{\prime})^{n+1} \rightarrow C/L^{n+1}$ is a thickening of discrete rings for any $n \in \mathbb{N}$.
    By the first paragraph of the lemma, we may find a lift $B \rightarrow C^{\prime} / L^{\prime}$ for the thickening of discrete rings $C^{\prime} / L^{\prime} \rightarrow C/L$.
    This gives the base case of the induction.
    Now assume we found a compatible lift $B \rightarrow C^{\prime} / (L^{\prime})^{n}$.
    We need to find a lift $B \rightarrow C^{\prime} / (L^{\prime})^{n+1}$ compatible with both $B \rightarrow C^{\prime} / (L^{\prime})^{n}$ and $B \rightarrow C/L^{n+1}$.
    Consider the induced morphism
    \begin{align}
        \label{eq:fibre_product_inductive_step}
        C^{\prime} / (L^{\prime})^{n+1} \rightarrow \left( C/L^{n+1} \right) \times_{C/L^{n}} \left( C^{\prime} / (L^{\prime})^{n} \right) =: D.
    \end{align}
    Since $C/L^{n} = \left( C^{\prime} / (L^{\prime})^{n} \right) \otimes_{C^{\prime} / (L^{\prime})^{n+1}} \left( C/L^{n+1} \right)$ and every morphism in the tensor product is a thickening of discrete rings, we may apply Lemma \ref{lem:chinese_remainder} (with $R = C^{\prime} / (L^{\prime})^{n+1}$) to conclude that the morphism in (\ref{eq:fibre_product_inductive_step}) is a thickening of discrete rings.
    Furthermore, by the universal property of fibre products of discrete rings, there is a induced compatible morphism $B \rightarrow D$.
    By the first paragraph of the lemma, we may find a lift $B \rightarrow C^{\prime} / (L^{\prime})^{n+1}$ for the thickening of discrete rings $C^{\prime} / (L^{\prime})^{n+1} \rightarrow D$.
    By construction, this lift is compatible with both morphisms $C^{\prime} / (L^{\prime})^{n+1} \rightarrow C / L^{n+1}$ and $C^{\prime} / (L^{\prime})^{n+1} \rightarrow C^{\prime} / (L^{\prime})^{n}$.
    This completes the inductive step and hence the proof of the lemma.
\end{proof}

\begin{lemma}
    \label{lem:formally_smooth_thickening_isomorphism_affine}
    Let $\varphi : A \rightarrow B$ be a thickening of adic rings which is formally smooth of formal finite presentation.
    Then $\varphi$ is an isomorphism.
\end{lemma}

\begin{proof}
    Apply Lemma \ref{lem:formally_smooth_lift_affine} to the solid diagram
    \begin{equation}
        \begin{tikzcd}
            \label{diag:formally_smooth_thickening_isomorphism_affine}
            A \arrow[d, "\varphi"] \arrow[r, "\mathds{1}_A"] & A \arrow[d, "\varphi"] \\
            B \arrow[r, "\mathds{1}_B"]  \arrow[ur, dashed] & B
        \end{tikzcd}
    \end{equation}
    to get existence of the dashed morphism, i.e. an inverse.
\end{proof}

Finally, we show that the property of being a thickening descends along thickenings.

\begin{lemma}
    \label{lem:closed_immersion_descent_affine}
    Let
    \begin{equation}
        \label{diag:closed_immersion_descent_affine}
        \begin{tikzcd}
            A \arrow[d, "\varphi"] \arrow[r, "\psi"] & A_0 \arrow[d, "\varphi_0"] \\
            B \arrow[r, "\psi_0"] & B_0
        \end{tikzcd}
    \end{equation}
    be a co-Cartesian diagram of adic rings and assume that $\psi$ is a thickening of adic rings of formal finite presentation.
    Then if $\varphi_0$ is a surjective adic morphism, so is $\varphi$, and if $\varphi_0$ is a thickening, so is $\varphi$.
\end{lemma}

\begin{proof}
    The proof consists in a reduction to the case of discrete rings. \par
    
    We first show that $\varphi$ is an adic morphism.
    Let $I \subseteq A$ and $J \subseteq B$ be ideals of definition.
    By assumption, $\psi$ and $\varphi_0$ are adic morphisms, hence $I \cdot A_0$ and $I \cdot B_0$ are ideals of definition of $A_0$ and $B_0$ respectively.
    By the $(1) \rightarrow (3)$ implication of Lemma \ref{lem:characterisation_thickenings}, $\mathrm{rad}_{B} \, J = \psi_0^{-1} \left( \mathrm{rad}_{B_0} \, (I \cdot B_0) \right)$.
    We will show that
    \begin{align}
        \label{eq:adic_ideal_thickening_descent}
        \mathrm{rad}_{B} \, (I \cdot B) = \psi_0^{-1} \left( \mathrm{rad}_{B_0} \, (I \cdot B_0) \right).
    \end{align}
    If this holds then, since $I \cdot B$ and $J$ are finitely generated, the equality $\mathrm{rad}_{B} \, J = \mathrm{rad}_{B} \, (I \cdot B)$ shows that $I \cdot B$ is an ideal of definition of $B$. \par

    We now show (\ref{eq:adic_ideal_thickening_descent}).
    Let $K = \ker (A \rightarrow A_0)$ and let $L = \ker (B \rightarrow B_0)$.
    Since Diagram (\ref{lem:closed_immersion_descent_affine}) is co-Cartesian, $L = K \cdot B$.
    By the $(1) \rightarrow (2)$ implication of Lemma \ref{lem:characterisation_thickenings}, $K \subseteq \mathrm{rad}_A \, I$.
    It follows that
    \begin{align}
        \label{eq:kernel_radical_inclusion}
        L = K \cdot B \subseteq \left( \mathrm{rad}_A \, I \right) \cdot B \subseteq \mathrm{rad}_B \, (I \cdot B).
    \end{align}
    By construction, $\psi_0 ( I \cdot B) = I \cdot B_0$, hence
    \begin{align}
        \label{eq:radical_equality}
        \psi_0^{-1} \left( \mathrm{rad}_{B_0} \, (I \cdot B_0) \right) &= \mathrm{rad}_{B} \, \left( \psi_0^{-1} \left(\psi_0 \left( I \cdot B \right) \right) \right) \nonumber \\
        &= \mathrm{rad}_{B} \, \left( I \cdot B + L \right) \nonumber \\
        &= \mathrm{rad}_{B} \, \left( \mathrm{rad}_{B} (I \cdot B) + \mathrm{rad}_{B} (L) \right) \nonumber \\
        &= \mathrm{rad}_{B} \, (I \cdot B),
    \end{align}
    where the last equality is a consequence of (\ref{eq:kernel_radical_inclusion}).

    Now we show that $\varphi$ is surjective.
    Since we already know that $\varphi$ is adic, it suffices to show that the induced morphism $\varphi_0 : A / I \rightarrow B / (I \cdot B)$ is surjective (part (1) of \cite[\href{https://stacks.math.columbia.edu/tag/0315}{Lemma 0315}]{stacks-project}).
    To this end, we may base change Diagram (\ref{diag:closed_immersion_descent_affine}) by the morphism $A \rightarrow A / I$ and obtain morphisms of discrete rings.
    Surjectivity of $\varphi_0$ now follows from the assumptions and \cite[\href{https://stacks.math.columbia.edu/tag/0896}{Lemma 0896}]{stacks-project} for discrete rings. \par

    Clearly, if $\varphi_0$ induces a homeomorphism of topological spaces, so does $\varphi$.
\end{proof}

\subsection{Morphisms of formal schemes}
\label{subsec:presentation_global}

\begin{definition}
    \label{def:finite_presentation}
	A morphism of formal schemes $f : X \rightarrow Y$ is \emph{locally of formal finite presentation} if for any affine open set $V = \spf A \subseteq Y$ and $U = \spf B \subseteq f^{-1}(V) \subseteq X$, the induced morphism of adic rings $A \rightarrow B$ is of formal finite presentation.
\end{definition}

Being locally of formal finite presentation is an affine-local property.

\begin{lemma}
    \label{lem:finite_presentation_local}
	Suppose that $f : X \rightarrow Y$ is a morphism of formal schemes and assume there exists an affine open covering $\{V_i = \spf A_i\}$ of $Y$ and an affine open covering $\{U_{i_j} = \spf B_{i_j}\}$ of $X$ such that $U_{i_j} \subseteq f^{-1}(V_i)$ and the induced morphism of adic rings $A_i \rightarrow B_{i_j}$ is of formal finite presentation for all $i$ and $i_j$.
	Then $f$ is locally of formal finite presentation.
\end{lemma}

\begin{proof}
	This follows immediately from Lemma \ref{lem:finite_presentation_affine_local} (see \cite[\href{https://stacks.math.columbia.edu/tag/01ST}{Lemma 01ST}]{stacks-project}) after employing the fact that $X$ and $Y$ have a basis of affine open formal subschemes.
\end{proof}

\begin{lemma}
    \label{lem:base_change_finite_presentation}
	The base change of a morphism of formal schemes locally of formal finite presentation is locally of formal finite presentation.
\end{lemma}

\begin{proof}
    Use Lemma \ref{lem:finite_presentation_local} to reduce to the affine case and then apply Lemma \ref{lem:base_change_finite_presentation_affine}.
\end{proof}

It is clear that the composition of morphisms of formal schemes locally of formal finite presentation is locally of formal finite presentation.

\begin{lemma}
    \label{lem:finite_presentation_factor}
    Let
    \begin{align}
        \label{eq:composition_finite_presentation}
        X \xrightarrow{f} Y \xrightarrow{g} Z
    \end{align}
    be morphisms of formal schemes.
    If $g$ and $g \circ f$ are locally of formal finite presentation, so is $f$.
\end{lemma}

\begin{proof}
    Use Lemma \ref{lem:finite_presentation_local} to reduce to the affine case and then apply Lemma \ref{lem:finite_presentation_factor_affine}.
\end{proof}

\begin{lemma}
    \label{lem:noetherian_ascent}
    Let $f : X \rightarrow Y$ be a morphism of formal schemes locally of formal finite presentation.
    Then if $Y$ is locally Noetherian, so is $X$.
\end{lemma}

\begin{proof}
    Being locally Noetherian can be checked on an affine open covering of $X$.
    Now use Lemma \ref{lem:finite_presentation_local} to reduce to the affine case and then apply Lemma \ref{lem:noetherian_ascent_affine}.
\end{proof}

We conclude by globalising the auxiliary results of the previous subsection, where possible.

\begin{lemma}
    \label{lem:formally_etale_lift}
    Let $f : X \rightarrow Y$ be a formally étale morphism of formal schemes locally of formal finite presentation and let $T \rightarrow T^{\prime}$ be a thickening of formal schemes fitting in a commutative solid diagram
    \begin{equation}
        \label{diag:formally_etale_lift}
        \begin{tikzcd}
            T \arrow[d] \arrow[r] & X \arrow[d, "f"] \\
            T^{\prime} \arrow[r]  \arrow[ur, dashed] & Y.
        \end{tikzcd}
    \end{equation}
    Then there exists a unique dashed morphism making the diagram commute.
\end{lemma}

\begin{proof}
    Since $f$ is formally étale, if a lift $T^{\prime} \rightarrow X$ exists, it is unique.
    Thus, it suffices to construct a lift affine locally and uniqueness will ensure local lifts may be glued together (see proof of \cite[\href{https://stacks.math.columbia.edu/tag/04FD}{Lemma 04FD}]{stacks-project}).
    Since affine open formal subschemes form a basis for the topology of formal schemes, we may assume that all formal schemes are affine.
    In this case, a lift $T^{\prime} \rightarrow X$ exists by Lemma \ref{lem:formally_smooth_lift_affine}.
\end{proof}

\begin{lemma}
    \label{lem:formally_smooth_thickening_isomorphism}
    Let $f : X \rightarrow Y$ be a formally smooth thickening of formal schemes locally of formal finite presentation.
    Then $f$ is an isomorphism.
\end{lemma}

\begin{proof}
    Being an isomorphism is local on the target hence the claim follows from Lemma \ref{lem:formally_smooth_thickening_isomorphism_affine}.
\end{proof}

\begin{lemma}
    \label{lem:closed_immersion_descent}
    Let
    \begin{equation}
        \label{diag:closed_immersion_descent}
        \begin{tikzcd}
            X_0 \arrow[d, "f_0"] \arrow[r, "g_0"] & X \arrow[d, "f"] \\
            Y_0 \arrow[r, "g"] & Y
        \end{tikzcd}
    \end{equation}
    be a Cartesian diagram of formal schemes and assume that $g$ is a thickening of formal schemes locally of formal finite presentation.
    Then, if $f_0$ is a closed immersion, so is $f$, and if $f_0$ is a thickening, so is $f$.
\end{lemma}

\begin{proof}
    We reduce to the affine case.
    Being a closed immersion or a thickening can be checked locally on the target (Lemma \ref{lem:closed_immersion_local} and Lemma \ref{lem:thickening_local}), hence suppose $Y$ is affine.
    If this is the case, $Y_0$ and $X_0$ are also affine.
    Furthermore, if $g$ is a thickening, so is $g_0$ (Lemma \ref{lem:base_change_thickenings}), hence $X$ is also affine (Lemma \ref{lem:image_thickening}).
    Now the result follows from Lemma \ref{lem:closed_immersion_descent_affine}.
\end{proof}

%% file: files/section5.tex
\section{Infinitesimal neighbourhoods}
\label{sec:infinitesimal_neighbourhoods}

\subsection{Affine infinitesimal neighbourhoods}
\label{subsec:neighbourhoods_local}

We now give the local definition of infinitesimal neighbourhoods and provide a local construction for surjective adic morphisms.

\begin{definition}
    Let $\varphi : A \rightarrow B$ be a morphism of adic rings,
    The \emph{affine infinitesimal neighbourhood} of $\varphi$ is a thickening of adic rings $\hat{\varphi} : \hat{A}_{\varphi} \rightarrow B$, or simply $\hat{A}$, factorising $\varphi$ with the following universal property: for all thickenings of adic rings $\varphi^{\prime} : C^{\prime} \rightarrow C$ and for all morphisms of adic rings $\psi : A \rightarrow C^{\prime}$ and $\psi^{\prime} : B \rightarrow C$ such that the solid diagram
    \begin{equation}
        \label{diag:affine_infinitesimal_neighbourhood}
        \begin{tikzcd}
            A \arrow[d, "\psi"] \arrow[r] \arrow[rr, "\varphi", bend left]& \hat{A}_{\varphi} \arrow[r, "\hat{\varphi}"] \arrow[dl, dashed] & B \arrow[d, "\psi^{\prime}"] \\
            C^{\prime} \arrow[rr, "\varphi^{\prime}"] & & C
        \end{tikzcd}
    \end{equation}
    is commutative, there exists a unique dashed morphism preserving commutativity of the diagram.
    By general category theory, if the affine infinitesimal neighbourhood exists, it is unique.
\end{definition}

\begin{proposition}
    \label{prop:existence_affine_infinitesimal_neighbourhood}
	Let $\varphi : A \rightarrow B$ be a surjective adic morphism of formal finite presentation.
	Let $I$ be an ideal of definition of $A$ and let $K = \ker \varphi$.
	Define
	\begin{align}
	    \label{eq:def_affine_infinitesimal_neighbourhood}
		\hat{\varphi} : \hat{A}_{\varphi} := \lim{n \in \mathbb{N}} A / \left( I+K \right)^{n+1} \rightarrow \lim{n \in \mathbb{N}} B / \left( I \cdot B \right)^{n+1} = B.
	\end{align}
	Then $\hat{\varphi}$ is the affine infinitesimal neighbourhood of $\varphi$.
\end{proposition}

\begin{proof}
	We first show that $\hat{\varphi}$ is a thickening of adic rings.
	Note that $\hat{A}_{\varphi}$ is the $(I + K)$-adic completion of $A$.
	Since $\varphi$ is adic and of formal finite presentation, the morphism $A \rightarrow B / \left(I \cdot B \right)$ is of finite presentation.
	As a result, its kernel, that is $(I+K)$, is a finitely generated ideal of $A$.
	This yields that $\hat{A}_{\varphi}$ is an adic ring and $(I+K) \cdot \hat{A}_{\varphi}$ is a finitely generated ideal of definition (part (1) of \cite[\href{https://stacks.math.columbia.edu/tag/05GG}{Lemma 05GG}]{stacks-project}).
	By construction, $\hat{\varphi}$ is surjective and adic and there exists a factorisation $\varphi : A \rightarrow \hat{A}_{\varphi} \rightarrow B$.
	In order to show $\hat{\varphi}$ is a thickening, it suffices to show $\hat{\varphi}$ is a homeomorphism.
    To this end, observe that, in the following commutative diagram
    \begin{equation}
        \label{diag:neighbourhood_homeomorphism}
        \begin{tikzcd}
            \hat{A}_{\varphi} \arrow[d] \arrow[r, "\hat{\varphi}"] & B \arrow[d] \\
            A / \left( I + K \right) \arrow[r, "\sim"] & B / \left( I \cdot B \right),
        \end{tikzcd}
    \end{equation}
    the vertical morphisms induce homeomorphisms on formal spectra and the lower horizontal morphism is an isomorphism. \par
	
	We now show the universal property.
    Let $L^{\prime} \subseteq C^{\prime}$ and $L \subseteq C$ be ideals of definition.
	In order to show existence and uniqueness of the dashed morphism, by the universal property of completion, it is necessary and sufficient to show that the morphism $\psi : A \rightarrow C^{\prime}$ is continuous when $A$ is endowed with the $(I + K)$-adic topology.
    Since $(I + K)$ is finitely generated, this is equivalent to requiring
    \begin{equation}
        \label{eq:continuous_universal_property_i}
        I + K \subseteq \mathrm{rad}_A \, \left( \psi^{-1} \left( L^{\prime} \right) \right).
    \end{equation}
    By assumption $\varphi^{\prime} : C^{\prime} \rightarrow C$ is a thickening, hence the implication $(1) \rightarrow (3)$ of Lemma \ref{lem:characterisation_thickenings} implies that
    \begin{align}
        \label{eq:continuous_universal_property_ii}
        \mathrm{rad}_{C^{\prime}} \, L^{\prime} = \left( \varphi^{\prime} \right)^{-1} \left( \mathrm{rad}_C \, L \right).
    \end{align}
    Now, using (\ref{eq:continuous_universal_property_ii}) and the fact that $\psi^{\prime}$ is continuous, it follows that
	\begin{align}
	    \label{eq:continuous_universal_property_iii}
		\mathrm{rad}_A \, \left( \psi^{-1} \left( L^{\prime} \right) \right) &= \psi^{-1} \left( \mathrm{rad}_{C^{\prime}} \, L^{\prime} \right) \nonumber \\
		&= \left( \varphi^{\prime} \circ \psi \right)^{-1} \left( \mathrm{rad}_{C} \, L \right) \nonumber \\
		&= \left( \psi^{\prime} \circ \varphi \right)^{-1} \left( \mathrm{rad}_{C} \, L \right) \nonumber \\
		&\supseteq \varphi^{-1} \left( I \cdot B \right) = (I + K)
	\end{align}
	This shows that the universal property holds.
\end{proof}

\begin{lemma}
    \label{lem:noetherian_neighbourhood_affine}
	Let $\varphi : A \rightarrow B$ be a surjective adic morphism of formal finite presentation.
    Then $\hat{\varphi} : \hat{A} \rightarrow B$ and $A \rightarrow \hat{A}$ are of formal finite presentation, and if $B$ is Noetherian, so is $\hat{A}$.
\end{lemma}

\begin{proof}
    In the notation of the proof of Proposition \ref{prop:existence_affine_infinitesimal_neighbourhood}, $\hat{A}_{\varphi}$ is the $(I+K)$-adic completion of $A$.
    Since $(I+K)$ is a finitely generated ideal of $A$, the morphism $A \rightarrow \hat{A}$ is of formal finite presentation.
    By Lemma \ref{lem:finite_presentation_factor_affine}, $\hat{\varphi}$ is also of formal finite presentation.
    Finally, if $A/(I+K) = B$ is Noetherian, so is $\hat{A}$ (\cite[\href{https://stacks.math.columbia.edu/tag/05GH}{Lemma 05GH}]{stacks-project}).
\end{proof}

\subsection{Infinitesimal neighbourhoods}
\label{subsec:neighbourhoods_global}

Now we finally define infinitesimal neighbourhoods in full generality and show that they exist for formally unramified adic morphisms locally of formal finite presentation.

\begin{definition}
    \label{def:infinitesimal_neighbourhood}
    Let $f : X \rightarrow Y$ be a morphism of formal schemes.
    The \emph{infinitesimal neighbourhood} of $f$ is a thickening of formal schemes $\hat{f} : X \rightarrow \hat{Y}$ factorising $f$ with the following universal property: for all thickenings of formal schemes $T \rightarrow T^{\prime}$ and for all morphisms $T \rightarrow X$ and $T^{\prime} \rightarrow Y$ such that the solid diagram
    \begin{equation}
        \label{diag:infinitesimal_neighbourhood}
        \begin{tikzcd}
            T \arrow[d] \arrow[rr] & & T^{\prime} \arrow[d] \arrow[dl, dashed] \\
            X \arrow[r, "\hat{f}"] & \hat{Y} \arrow[r, "\iota"] & Y 
        \end{tikzcd}
    \end{equation}
    is commutative, there exists a unique dashed morphism preserving commutativity of the diagram.
    By general category theory, if the infinitesimal neighbourhood exists, it is unique.
\end{definition}

Firstly, we show that the affine infinitesimal neighbourhood of a morphism of affine formal schemes is the infinitesimal neighbourhood in the category of formal schemes.

\begin{lemma}
    \label{lem:affine_infinitesimal_neighbourhood_formal_schemes}
    Let $f: X \rightarrow Y$ be a morphism of affine formal schemes and assume that the affine infinitesimal neighbourhood
    \begin{align}
        \label{eq:affine_infinitesimal_neighbourhood}
        X \xrightarrow{\hat{f}} \hat{Y} \xrightarrow{\iota} Y
    \end{align}
    of $f$ exists in the category of affine formal schemes.
    Then $\hat{f}$ is the infinitesimal neighbourhood of $f$ in the category of formal schemes.
\end{lemma}

\begin{proof}
    We need to show that the dashed morphism in Diagram (\ref{diag:infinitesimal_neighbourhood}) exists and is unique.
    Suppose firstly that $T^{\prime}$ is an affine formal scheme.
    Since $T \rightarrow T^{\prime}$ is a thickening, $T$ is also an affine formal scheme.
    Therefore all the schemes are affine and the dashed morphism exists by definition of affine infinitesimal neighbourhood. \par
    
    Next cover $T^{\prime}$ by affine open formal subschemes $T^{\prime}_i$, for $i \in I$.
    By the previous paragraph there exists a unique morphism $T^{\prime}_i \rightarrow \hat{Y}$ for all $i \in I$.
    Over an affine open formal subscheme of the pairwise intersection $T^{\prime}_i \cap T^{\prime}_j$, for $i$ and $j$ in $I$, the restrictions of $T^{\prime}_i \rightarrow \hat{Y}$ and $T^{\prime}_j \rightarrow \hat{Y}$ agree by uniqueness.
    Therefore the morphisms glue to give a unique morphism $T^{\prime} \rightarrow \hat{Y}$.
\end{proof}

Next, we show that infinitesimal neighbourhoods are not affected by formally étale morphisms.

\begin{lemma}
    \label{lem:infinitesimal_neighbourhood_formally_etale}
    Let
    \begin{align}
        X \xrightarrow{f} Y \xrightarrow{g} Z
    \end{align}
    be morphisms of formal schemes and suppose that $g$ is formally étale and locally of formal finite presentation.
    Then, if the infinitesimal neighbourhood of $f$ exists, so does the infinitesimal neighbourhood of $g \circ f$, and is equal to the infinitesimal neighbourhood of $f$.
\end{lemma}

\begin{proof}
    Consider the solid commutative diagram
    \begin{equation}
        \label{diag:infinitesimal_neighbourhood_formally_etale}
        \begin{tikzcd}
            T \arrow[d] \arrow[rrr] & & & T^{\prime} \arrow[d] \arrow[dl, dotted] \arrow[dll, dashed] \\
            X \arrow[r, "\hat{f}"] & \hat{Y} \arrow[r, "\iota"] & Y \arrow[r, "g"] & Z. 
        \end{tikzcd}
    \end{equation}
    We show that $\hat{f}$ satisfies the universal property of the infinitesimal neighbourhood of $g \circ f$.
    By Lemma \ref{lem:formally_etale_lift}, there exists a unique compatible dotted morphism.
    By the universal property of the infinitesimal neighbourhood of $f$, there exists a unique compatible dashed morphism.
\end{proof}

Next, we show that infinitesimal neighbourhoods commute with fibre products, when they exist.
This is not surprising as they are both categorical limits.

\begin{lemma}
    \label{lem:infinitesimal_neighbourhood_fibre_products}
    Let
    \begin{equation}
        \begin{tikzcd}
            X^{\prime} \arrow[d] \arrow[r, "f^{\prime}"] & Y^{\prime} \arrow[d]\\
            X \arrow[r, "f"] & Y
        \end{tikzcd}
    \end{equation}
    be a Cartesian diagram of formal schemes.
    Then, if the infinitesimal neighbourhood $\hat{Y}$ of $f$ exists, so does the infinitesimal neighbourhood of $f^{\prime}$, and is given by $\hat{Y} \times_{Y} Y^{\prime}$.
\end{lemma}

\begin{proof}
    Define $\hat{Y}^{\prime} := \hat{Y} \times_{Y} Y^{\prime}$ and consider the solid commutative diagram of formal schemes
    \begin{equation}
        \label{diag:infinitesimal_neighbourhood_fibre_products}
        \begin{tikzcd}
            T \arrow[d] \arrow[rr] & & T^{\prime} \arrow[d] \arrow[dl, dashed] \\
            X^{\prime} \arrow[d] \arrow[r] & \hat{Y}^{\prime} \arrow[d] \arrow[r] & Y^{\prime} \arrow[d] \\
            X \arrow[r] & \hat{Y} \arrow[r] & Y,
        \end{tikzcd}
    \end{equation}
    where the two bottom squares are Cartesian and $T \rightarrow T^{\prime}$ is a thickening of formal schemes.
    Suppose that $\hat{Y}$ is the infinitesimal neighbourhood of $f$, we show that the fibre product $\hat{Y}^{\prime}$ is the infinitesimal neighbourhood of $f^{\prime}$.
    To this end, it suffices to show that a compatible dashed morphism exists and is unique.
    By the universal property of the infinitesimal neighbourhood $\hat{Y}$, there exists a unique compatible morphism $T^{\prime} \rightarrow \hat{Y}$.
    By the universal property of the fibre product $\hat{Y}^{\prime}$, there exists a unique compatible dashed morphism $T^{\prime} \rightarrow \hat{Y}^{\prime}$.
\end{proof}

Next, we show that it suffices to construct infinitesimal neighbourhoods on an open covering.

\begin{lemma}
    \label{lem:infinitesimal_neighbourhood_open_covering}
    Let $f : X \rightarrow Y$ be a morphism of formal schemes.
    Assume that, for all $x \in X$, there exist open formal subschemes $V \subseteq Y$ and $U \subseteq X$, with $x \in U \subseteq f^{-1}(V)$, such that the infinitesimal neighbourhood of the restriction $f_U : U \rightarrow V$ exists.
    Then the infinitesimal neighbourhood of $f$ exists.
\end{lemma}

\begin{proof}
    Let $\{ X_i \}_{i \in I}$ and $\{ Y_i \}_{i \in I}$ denote the collections of such open formal subschemes, where $I$ indexes the points of $X$, and let $f_i : X_i \rightarrow Y_i$ denote the restriction of $f$.
    By assumption, the infinitesimal neighbourhood $\hat{Y}_i$ of $f_i$ exists. \par
    
    We glue the collection of affine formal schemes $\{\hat{Y}_i\}_{i \in I}$.
    Let $Y_{ijk}$ denote the triple intersection of $Y_i$, $Y_j$ and $Y_k$ for $i,j,k \in I$.
    If the indices are repeated, the double intersection $Y_{ijj}$ is simply written as $Y_{ij}$.
    Let $f_{ijk} : X_{ijk} \rightarrow Y_{ijk}$ be the base change of $f_i$ by the open immersion $Y_{ijk} \rightarrow Y_i$.
    By Lemma \ref{lem:infinitesimal_neighbourhood_fibre_products}, the infinitesimal neighbourhood $\hat{Y}_{ijk}$ of $f_{ijk}$ exists and is an open formal subscheme of $\hat{Y}_i$.
    Using the universal property of infinitesimal neighbourhoods, it is clear that $\hat{Y}_{ij}$ and $\hat{Y}_{ji}$ are uniquely isomorphic, both being infinitesimal neighbourhoods of $f_{ij} = f_{ji}$.
    By a similar argument, there exist unique canonical isomorphisms
    \begin{align}
        \hat{Y}_{ijk} \cong \hat{Y}_{kij} \cong \hat{Y}_{jki}.
    \end{align}
    As a result, the affine formal schemes $\{ \hat{Y}_i \}_{i \in I}$ can be glued to obtain a formal scheme $\hat{Y}$ and a factorisation
    \begin{align}
        f : X \xrightarrow{\hat{f}} \hat{Y} \xrightarrow{\iota} Y.
    \end{align}
    By construction and by Lemma \ref{lem:thickening_local}, it follows that $\hat{f}$ is a thickening of formal schemes. \par
    
    Finally, we verify that $\hat{Y}$ satisfies the universal property of infinitesimal neighbourhoods.
    For any $i,j \in I$, let $T^{\prime}_{ij} \rightarrow Y_{ij}$ be the the base change of $T^{\prime} \rightarrow Y$ by the open immersion $Y_{ij} \rightarrow Y$.
    By the universal property of $\hat{Y}_{ij}$, there exists a unique compatible morphism $T^{\prime}_{ij} \rightarrow \hat{Y}_{ij}$.
    Therefore, the morphisms $T^{\prime}_{i} \rightarrow \hat{Y}$ glue to give a unique morphism $T^{\prime} \rightarrow \hat{Y}$.
\end{proof}

Next, we show auxiliary properties of the morphism $\iota$ of Definition \ref{def:infinitesimal_neighbourhood}.

\begin{lemma}
    \label{lem:iota_properties}
    Let $f : X \rightarrow Y$ be a morphism of formal schemes and suppose that the infinitesimal neighbourhood 
    \begin{align}
        \label{eq:infinitesimal_neighbourhood}
        X \xrightarrow{\hat{f}} \hat{Y} \xrightarrow{\iota} Y
    \end{align}
    of $f$ exists.
    Then
    \begin{enumerate}
        \item $\iota$ is formally étale;
        \item if $f$ is a monomorphism, so is $\iota$; and
        \item if $f$ is affine, so is $\iota$.
    \end{enumerate}
\end{lemma}

\begin{proof}
    \begin{enumerate}[itemsep=1em]
        \item Consider the commutative diagram of solid arrows
        \begin{equation}
            \label{diag:iota_formally_smooth}
            \begin{tikzcd}
                X \times_{\hat{Y}} T \arrow[d] \arrow[r] & T \arrow[d] \arrow[r] & T^{\prime} \arrow[d] \arrow[dl, dashed] \\
                X \arrow[r, "\hat{f}"] & \hat{Y} \arrow[r, "\iota"] & Y 
            \end{tikzcd}
        \end{equation}
        with $T \rightarrow T^{\prime}$ a first order thickening of affine schemes.
        We have to show existence and uniqueness of the dashed morphism.
        Note that, since $\hat{f}$ and $T \rightarrow T^{\prime}$ are thickenings, the composition $X \times_{\hat{Y}} T \rightarrow T^{\prime}$ is a thickening.
        Now existence and uniqueness of the dashed morphism follow from the universal property of infinitesimal neighbourhoods.

        \item Let $T \rightrightarrows \hat{Y}$ be two morphisms which are equal after composing with $\iota$.
        Since $\hat{Y}_{\mathrm{red}} = X_{\mathrm{red}}$, we obtain two morphisms $T_{\mathrm{red}} \rightrightarrows X$ which are equal after composing with $f$.
        Since $f$ is a monomorphism, the two morphisms $T_{\mathrm{red}} \rightrightarrows X$ are equal.
        Hence consider the commutative diagram
        \begin{equation}
            \label{diag:iota_monomorphism}
            \begin{tikzcd}
                T_{\mathrm{red}} \arrow[d] \arrow[rr] & & T \arrow[d] \arrow[dl, shift left = 0.5ex] \arrow[dl, shift right = 0.5ex] \\
                X \arrow[r, "\hat{f}"] & \hat{Y} \arrow[r, "\iota"] & Y.
            \end{tikzcd}
        \end{equation}        
        Since $T_{\mathrm{red}} \rightarrow T$ is a thickening, there is at most one morphism $T \rightarrow \hat{Y}$.

        \item By assumption, if $Y$ is affine, so is $X$.
        Now Lemma \ref{lem:affine_infinitesimal_neighbourhood_formal_schemes} implies that $\hat{Y}$ is also affine.
    \end{enumerate}
    This concludes the proof of the lemma.
\end{proof}

We finally prove existence of infinitesimal neighbourhoods and the Main Theorem.

\begin{proposition}
    \label{prop:infinitesimal_neighbourhoods}
    Let $f : X \rightarrow Y$ be a formally unramified adic morphism of formal schemes locally of formal finite presentation.
    Then the infinitesimal neighbourhood of $f$ exists.
\end{proposition}

\begin{proof}
    By Lemma \ref{lem:infinitesimal_neighbourhood_open_covering}, it suffices to construct the infinitesimal neighbourhood of $f$ locally.
    Hence we may assume that $X = \spf B$ and $Y = \spf A$ are affine formal schemes, and $f$ is induced by $\varphi : A \rightarrow B$.
    This morphism is formally unramified, adic and of formal finite presentation. \par

    Let $I$ and $J$ be ideals of definition of $A$ and $B$ respectively.
    By assumption and the fact that $B \rightarrow B/J$ is a surjective adic morphism, we have that $A \rightarrow B/J$ is formally unramified and of finite presentation, hence unramified (\cite[\href{https://stacks.math.columbia.edu/tag/02H7}{Section 02H7}]{stacks-project}).
    By the étale local structure theorem of unramified morphisms, for any $x \in \spec B/J$, we may find an element $s \in B/J$, such that $x \in \spec (B/J)_{s}$, and a factorisation
    \begin{align}
        \label{eq:formally_unramified_factorisation_discrete}
        A \xrightarrow{\psi} E \xrightarrow{\varphi^{\prime}} (B/J)_{s}
    \end{align}
    of $\varphi$, such that $\psi$ is an étale morphism (as a morphism of discrete rings) and $\varphi^{\prime}$ is surjective (\cite[\href{https://stacks.math.columbia.edu/tag/0395}{Proposition 0395}]{stacks-project}).
    Furthermore, $\varphi^{\prime}$ is of formal finite presentation (Lemma \ref{lem:finite_presentation_factor_affine}). Let $K$ be its kernel. \par

    Let $\hat{\varphi}^{\prime} : \hat{E} \rightarrow (B/J)_{s}$ be the infinitesimal neighbourhood of $\varphi^{\prime}$ (Proposition \ref{prop:existence_affine_infinitesimal_neighbourhood}).
    This is a thickening of adic rings of formal finite presentation (Lemma \ref{lem:noetherian_neighbourhood_affine}).
    We now show that the composition $A \rightarrow E \rightarrow \hat{E}$ is a morphism of adic rings which is formally étale and of formal finite presentation.
    By (\ref{eq:def_affine_infinitesimal_neighbourhood}), we know that $\hat{E}$ is defined as the $K$-adic completion of $E$.
    By construction, $\psi^{-1}(K) = \varphi^{-1}(J) \subseteq A$.
    Since $\varphi$ is continuous, it follows that the morphism $A \rightarrow \hat{E}$ is continuous.
    Furthermore, since $A \rightarrow E$ is étale, the morphism $A \rightarrow \hat{E}$ is formally étale (part (1) of Lemma \ref{lem:iota_properties}) and of formal finite presentation (Lemma \ref{lem:noetherian_neighbourhood_affine}). \par

    Now consider the morphism $B \rightarrow B/J$ and let $\tilde{s} \in B$ be an element mapping to $s$.
    We have an induced morphism $B_{\{\tilde{s}\}} \rightarrow (B/J)_{s}$.
    This is the base change of $B \rightarrow B/J$ by the localisation morphism $B \rightarrow B_{\{\tilde{s}\}}$.
    Furthermore, $x \in \spf B_{\{\tilde{s}\}}$, thus, up to replacing $X$, we may assume that $\tilde{s} = 1$.
    We obtain a solid commutative diagram
    \begin{equation}
        \label{diag:formally_unramified_factorisation}
        \begin{tikzcd}
            A \arrow[d] \arrow[r, "\varphi"] & B \arrow[d] \\
            \hat{E} \arrow[r, "\hat{\varphi}^{\prime}"] \arrow[ur, dashed, "\varphi^{\prime\prime}" ] & B/J,
        \end{tikzcd}
    \end{equation}
    where $A \rightarrow \hat{E}$ is formally étale of formal finite presentation and $B \rightarrow B/J$ is a thickening.
    By Lemma \ref{lem:formally_smooth_lift_affine}, there exists a compatible dashed morphism $\varphi^{\prime\prime} : \hat{E} \rightarrow B$.
    We note that, since $\varphi$ is adic and of formal finite presentation, the same conditions hold for $\varphi^{\prime\prime}$.
    The first condition follows directly from \cite[\href{https://stacks.math.columbia.edu/tag/0GXB}{Lemma 0GXB}]{stacks-project}, and the second condition from Lemma \ref{lem:finite_presentation_factor} and the fact that $A \rightarrow \hat{E}$ is of formal finite presentation.
    Now surjectivity of $\hat{\varphi}^{\prime}$ and part (1) of \cite[\href{https://stacks.math.columbia.edu/tag/0315}{Lemma 0315}]{stacks-project} imply that $\varphi^{\prime\prime}$ is surjective. \par

    By Lemma \ref{lem:infinitesimal_neighbourhood_formally_etale} and the fact that $A \rightarrow \hat{E}$ is a formally étale morphism of formal finite presentation, in order to construct the infinitesimal neighbourhood of $\varphi$, it is sufficient to construct the infinitesimal neighbourhood of $\varphi^{\prime\prime}$.
    Since $\varphi^{\prime\prime}$ is a surjective adic morphism of formal finite presentation, its affine infinitesimal neighbourhood exists by Proposition \ref{prop:existence_affine_infinitesimal_neighbourhood}, and it is the infinitesimal neighbourhood by Lemma \ref{lem:affine_infinitesimal_neighbourhood_formal_schemes}.
\end{proof}

\begin{lemma}
    \label{lem:noetherian_neighbourhood}
    Let $f : X \rightarrow Y$ be a formally unramified morphism of formal schemes locally of formal finite presentation.
    Then $\hat{f}$ and $\iota : \hat{Y} \rightarrow Y$ are locally of formal finite presentation, and if $X$ is locally Noetherian, so is $\hat{Y}$.
\end{lemma}

\begin{proof}
    We use the notation of the proof of Proposition \ref{prop:infinitesimal_neighbourhoods}.
    By Lemma \ref{lem:finite_presentation_factor}, in order to show that $\hat{f}$ is locally of formal finite presentation, it suffices to show that so is $\iota$.
    On affine open covering, $\iota$ is the formal spectrum of the composition of $A \rightarrow \hat{E}$ (Diagram (\ref{diag:formally_unramified_factorisation})) with the infinitesimal neighbourhood of $\hat{E} \rightarrow B$.
    The former morphism is of finite presentation, as remarked in the proof of Proposition \ref{prop:infinitesimal_neighbourhoods}, and the latter morphism is of finite presentation (Lemma \ref{lem:noetherian_neighbourhood_affine}).
    Since being locally of formal finite presentation is a local property (Lemma \ref{lem:finite_presentation_local}), we conclude that $\iota$ is locally of formal finite presentation.
    Finally, since being locally Noetherian can be checked on an open covering, Lemma \ref{lem:noetherian_neighbourhood_affine} and the local construction of $\hat{Y}$ imply that $\hat{Y}$ is locally Noetherian, if so is $X$.
\end{proof}

\begin{proof}
    [Proof of Main Theorem]
    By Proposition \ref{prop:infinitesimal_neighbourhoods}, the infinitesimal neighbourhood of $f$ exists.
    By Lemma \ref{lem:noetherian_neighbourhood}, $\hat{f}$ and $\iota$ are locally of formal finite presentation.
    By part (1) of Lemma \ref{lem:iota_properties}, $\iota$ is formally étale.
\end{proof}

\begin{example}
    \label{ex:immersions}
    Let $f : X \rightarrow Y$ be an immersion of formal schemes locally of formal finite presentation.
    Then it is formally unramified and adic, hence we may construct the infinitesimal neighbourhood $\hat{f} : X \rightarrow \hat{Y}$.
    This is locally given by completing with respect to an ideal.
    Indeed, since $f$ is an immersion, $f$ is locally given by the formal spectrum of a surjective adic morphism $\varphi : A \rightarrow B$.
    Let $K$ be its kernel.
    By Proposition \ref{prop:existence_affine_infinitesimal_neighbourhood}, $\hat{Y}$ is locally given by the formal spectrum of the completion of $A$ with respect to the ideal $(I + K)$, where $I$ is an ideal of definition of $A$.
\end{example}